\documentclass{amsart}
\usepackage{graphicx}

\newtheorem{theorem}{Theorem}[section]
\newtheorem{lemma}[theorem]{Lemma}
\newtheorem{proposition}[theorem]{Proposition}

\begin{document}

\title{Additivity of free genus of knots}

\author{Makoto Ozawa}
\address{Department of Mathematics, School of Education, Waseda University,
1-6-1 Nishiwaseda, Shinjuku-ku, Tokyo 169-8050, Japan}
\email{ozawa@mn.waseda.ac.jp}

\begin{abstract}
We show that free genus of knots is additive under connected sum.
\end{abstract}

\maketitle

\section{Introduction}
Let $K$ be a knot in the 3-sphere $S^3$.
A Seifert surface $F$ for $K$ in $S^3$ is said to be {\it free} if the fundamental group $\pi _1(S^3-F)$ is a free group.
We note that all knots bound free Seifert surfaces, e.g. canonical Seifert surfaces constructed by the Seifert's algorithm.
We define the {\it free genus} $g_f(K)$ of $K$ as the minimal genus over all free Seifert surfaces for $K$ (\cite{K}).

For usual genus, Schubert  (\cite[2.10 Proposition]{S}) proved that genus of knots is additive under connected sum.
In general, the genus of a knot is not equal to the free genus of it.
In fact, free genus may have arbitrarily high gaps with genus (\cite{M}, \cite{KK}).

In this paper, we show the following theorem.

\begin{theorem}
For two knots $K_1$, $K_2$ in $S^3$, $g_f(K_1)+g_f(K_2)=g_f(K_1\# K_2)$.
\end{theorem}

\section{Preliminaries}
We can deform a Seifert surface $F$ by an isotopy so that $F\cap N(K)=N(\partial F;F)$.
We denote the exterior $cl(S^3-N(K))$ by $E(K)$, and the exterior $cl(S^3-N(F))$ or $cl(E(K)-N(F))$ by $E(F)$.
We have the following proposition.

\begin{proposition}$($\cite[5.2]{H}, \cite[IV.15]{J}, \cite[Lemma 2.2.]{O}$)$
A Seifert surface $F$ is free if and only if $E(F)$ is a handlebody.
\end{proposition}

We have the following inequality.

\begin{proposition}
$g_f(K_1)+g_f(K_2)\ge g_f(K_1\# K_2)$.
\end{proposition}

\begin{proof}
Let $F_i\ (i=1,2)$ be a free Seifert surface of minimal genus for $K_i$.
We construct a Seifert surface $F$ for $K_1\# K_2$ as the boundary connected sum of $F_1$ and $F_2$ naturally.
Then $E(F)$ is obtained by a boundary connected sum of $E(F_1)$ and $E(F_2)$.
Therefore the exterior of $F$ is a handlebody, and $F$ is free.
Hence we have the desired inequality. 
\end{proof}

We can specify the {\it $+$-side} and {\it $-$-side} of a Seifert surface $F$
for a knot $K$ by the orientability of $F$.
We say that a compressing disk $D$ for $F$ is a {\it $+$-compressing disk} (resp. {\it $-$-compressing disk}) if the collar of its boundary lies on the $+$-side (resp. $-$-side) of $F$, and $F$ is called {\it $+$-compressible} (resp. {\it $-$-compressible}) if $F$ has a {\it $+$-compressing disk} (resp. {\it $-$-compressing disk}).
A Seifert surface is said to be {\it weakly reducible} if there exist a $+$-compressing disk $D^+$ and a $-$-compressing disk $D^-$ for $F$ such that $\partial D^+\cap \partial D^-=\emptyset$. Otherwise $F$ is {\it strongly irreducible}.
The Seifert surface $F$ is {\it reducible} if $\partial D^+=\partial D^-$. Otherwise $F$ is {\it irreducible}.

\begin{proposition}
A free Seifert surface of minimal genus is irreducible.
\end{proposition}

\begin{proof}
Suppose that $F$ is reducible.
Then there exist a $+$-compressing disk $D^+$ and a $-$-compressing disk $D^-$ for $F$ such that $\partial D^+=\partial D^-$.
By a compression of $F$ along $D^+$ (this is the same to a compression along $D^-$), we have a new Seifert surface $F'$.
Since $E(F')$ is homeomorphic to a component of the manifold which is obtained by cutting $E(F)$ along $D^+\cup D^-$, it is a handlebody.
Hence $F'$ is free, but it has a lower genus than $F$.
This contradicts the minimality of $F$. 
\end{proof}

For a free Seifert surface $F$ of minimal genus for $K_1\# K_2$ and a decomposing sphere $S$ for the connected sum of $K_1$ and $K_2$, we will show ultimately that $S$ can be deformed by an isotopy so that $S$ intersects $F$ in a single arc, and we have the equality in Theorem 1.
To do this, we divide the proof of Theorem 1 into two cases;
(1) $F$ is strongly irreducible,
(2) $F$ is weakly reducible.
The case (1) is treated in the next section and we consider the case (2) in Section 4.

\section{Proof of Theorem 1 (strongly irreducible case)}
If a free Seifert surface $F$ of minimal genus for $K_1\# K_2$ is incompressible, then an innermost loop argument shows that a decomposing sphere $S$ for $K_1\# K_2$ can be deformed by an isotopy so that $S$ intersects $F$ in a single arc, and by Proposition 3, we have the equality in Theorem 1.

So, hereafter we suppose that $F$ is compressible and that in this section, $F$ is strongly irreducible.
Without loss of generality, we may assume that there is a $+$-compressing disk for $F$.
Let $\mathcal{D^+}$ be a $+$-compressing disk system for $F$, and let $F'$ be a surface obtained by compressing $F$ along $\mathcal{D^+}$.
We can retake $\mathcal{D^+}$ so that $F'$ is connected since $E(F)$ is a handlebody.
Take $\mathcal{D^+}$ to be maximal with respect to above conditions.
We deform $F'$ by an isotopy so that $F'\cap F=K$.
Put $A=\partial N(K_1\# K_2)-IntN(F)$, and let $H$ be a closed surface which is  obtained by pushing $F\cup A\cup F'$ into the interior of $E(F')$.
Let $A_0$ be a vertical annulus connecting a core of $A$ and a core of the copy of $A$ in $H$.
Then $H$ bounds a handlebody $V$ in $E(F')$ since $V$ is obtained from $E(F)$ by cutting along $\mathcal{D^+}$.
The remainder $W=E(F')-IntV$ is a compression body since it is obtained from $N(\partial E(F'); E(F'))$ by adding 1-handles $N(\mathcal{D^+})$.

\begin{figure}[htbp]
	\begin{center}
		\includegraphics[trim=5mm 0mm 5mm -5mm, width=.4\linewidth]{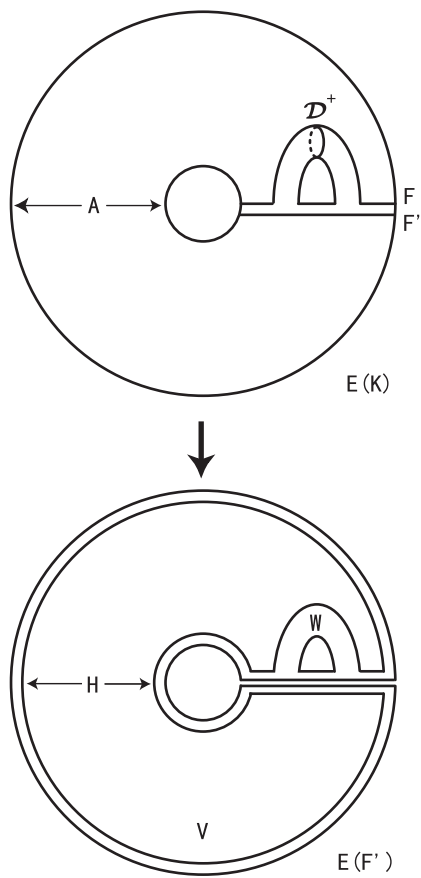}
	\end{center}
	\caption{Construction of a Heegaard splitting of $E(F')$}
\end{figure}

\begin{lemma}
$F'$ is incompressible in $S^3$.
\end{lemma}

\begin{proof}
We specify the $\pm$-side of $F'$ endowed from $F$ naturally.
Suppose that $F'$ is $+$-compressible, and let $E^+$ be a $+$-compressible disk for $F'$.
Then we can regard $E^+$ as a $\partial$-reducing disk for $E(F')$.
By applying our situation to \cite[Lemma 1.1.]{CG}, we may assume that $E^+\cap \mathcal{D^+}=\emptyset$.
If $\partial E^+$ separates $F'$, then $E^+$ cuts off a handlebody from $E(F')$, and there is a non-separating disk in it.
So, we may assume that $\partial E^+$ is non-separating in $F'$.
Then ${\mathcal{D^+}}\cup E^+$ is a $+$-compressing disk system satisfying the previous conditions.
This contradicts the maximality of $\mathcal{D^+}$.

Next, suppose that $F'$ is $-$-compressible, and let $E^-$ be a $-$-compressing disk for $F'$.
Then we can regard $E^-$ as a $\partial$-reducing disk for $E(F')$.
By applying our situation to \cite[Lemma 1.1.]{CG},
we may assume that $E^-\cap H=E^-\cap F$ is a single loop, and by exchanging $\mathcal{D^+}$ if it is necessary, that $E^-$ does not intersect $\mathcal{D^+}$.
But this contradicts the strongly irreducibility of $F$. 
\end{proof}

By Lemma 5,
We can deform the decomposing sphere $S$ by an isotopy so that $S$ intersects $F'$ in a single arc.
Put $E(S)=S\cap E(F')$.
Then $E(S)$ is a $\partial$-reducing disk for $E(F')$.
Otherwise, at least one of $K_1$ or $K_2$ is trivial, and Theorem 1 clearly holds.
By applying our situation to \cite[Lemma 1.1.]{CG},
we may assume that $E(S)$ intersects $H$ in a single loop, $E(S)$ intersects $A_0$ in two vertical arc, and by exchanging $\mathcal{D^+}$ under the previous conditions if it is necessary, that $E(S)$ does not intersects $\mathcal{D^+}$.
Then $S$ intersects $F$ in a single arc, hence we obtain the inequality $g_f(K_1)+g_f(K_2)\le g_f(K_1\# K_2)$.
This and Proposition 3 complete the Proof of Theorem 1.

\section{Proof of Theorem 1 (weakly reducible case)}
In this section, we consider the case that $F$ is weakly compressible.

We use the {\it Hayashi-Shimokawa (HS-) complexity} (\cite{HS}).
Here we review it.
Let $H$ be a closed (possibly disconnected) 2-manifold.
Put $w(H)=\{ genus(T)|T$ is a component of $H\}$, where this ``multi-set" may contain the same ordered pairs redundantly.
We order finite multi-sets as follows: arrange the elements of the multi-set in each one in monotonically non-increasing order, then compare the elements lexicographically.
We define the HS-complexity $c(H)$ as a multi-set obtained from $w(H)$ by deleting all the $0$ elements.
We order $c(H)$ in the same way as $w$.

Let $\alpha$ be a 1-submanifold of $H$.
Then let $\rho (H,\alpha )$ denote the closed 2-manifold obtained by cutting $H$ along $\alpha$ and capping off the resulting two boundary circles with disks.

Since $F$ is weakly reducible, there exist $+$-compressing disk $D^+$ and $-$-compressing disk $D^-$ for $F$ such that $\partial D^+\cap \partial D^-=\emptyset$.
If $c(\rho (F;\partial D^+\cup \partial D^-)) = c(\rho (F;\partial D^+))$, say, then $\partial D^-$ bounds a $+$-compressing disk for $F$.
Hence $F$ is reducible, and by Proposition 4, a contradiction.

Therefore there exist non-empty $+$-compressing disks system $\mathcal{D}^+$ and $-$-compressing disk system $\mathcal{D}^-$ for $F$ such that
\begin{itemize}
\item[(1)] $\partial \mathcal{D}^+\cap \partial \mathcal{D}^-=\emptyset$,
\item[(2)] $c(\rho (F;\partial \mathcal{D}^+\cup \partial \mathcal{D}^-))$ $<$ $c(\rho (F;\partial \mathcal{D}^+))$, $c(\rho (F;\partial \mathcal{D}^-))$,
\end{itemize}
and with $c(\rho (F;\partial \mathcal{D}^+\cup \partial \mathcal{D}^-))$ minimal subject to these conditions.
Moreover we take $\mathcal{D}^{\pm}$ so that $|\mathcal{D}^{\pm}|$ is minimal.

Let $F^{\pm}$ be a 2-manifold obtained by compressing $F$ along $\mathcal{D}^{\pm}$, and $F'$ be a 2-manifold obtained by compressing $F$ along $\mathcal{D}^+\cup \mathcal{D}^-$.
We deform $F^+$ and $F^-$ by an isotopy so that $F^+\cap F'\cap F^-=K$ and $F^{\pm}\cap N(K)=N(\partial F^{\pm}; F^{\pm})$.
Put $A=\partial N(K_1\# K_2)-IntN(F)$, and let $H$ be a closed 2-manifold which is obtained by pushing $F^+\cup A\cup F^-$ into the interior of $E(F')$.
Let $A_0$ be a vertical annulus connecting a core of $A$ and a core of the copy of $A$ in $H$.
Then $H$ bounds the union of handlebodies $V$ in $E(F')$ since $V$ is obtained from $E(F)$ by cutting along $\mathcal{D}^{\pm}$.
The remainder $W=E(F')-IntV$ is a union of compression bodies since it is obtained from $N(\partial E(F'); E(F'))$ by adding 1-handles $N(\mathcal{D}^{\pm})$.

\begin{figure}[htbp]
	\begin{center}
		\includegraphics[trim=5mm 0mm 5mm -5mm, width=.4\linewidth]{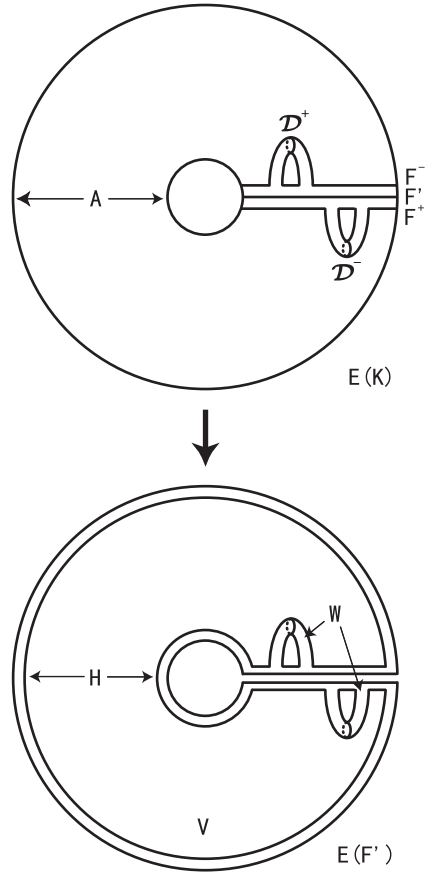}
	\end{center}
	\caption{Construction of a Heegaard splitting of $E(F')$}
\end{figure}

\begin{lemma}
There is no $2$-sphere component of $H$.
\end{lemma}

\begin{proof}
Suppose that there is a 2-sphere component $H_i$ of $H$.
We may assume that $H$ does not contain $A$, and there is a copy of some component of $\mathcal{D}^+$ in $H$.
Let ${\mathcal{D}}^+_s$ be a subsystem of $\mathcal{D}^+$ the union of whose boundaries separates $F$.
If there is no copy of $\mathcal{D}^-$ in $H_i$, then we delete any one of ${\mathcal{D}}^+_s$.
Then $\mathcal{D}^{\pm}$ holds the previous conditions, but this contradicts the minimality of $|{\mathcal{D}}^+|$.
If there is a copy of $\mathcal{D}^-$ in $H_i$, then there is a simple closed curve in $H_i$
which separates $N({\mathcal{D}^+}) \cap H_i$ from $N({\mathcal{D}^-}) \cap H_i$, and bounds a $+$-compressing disk and $-$-compressing disk for $F$.
Hence $F$ is reducible, but this contradicts Proposition 4.

\end{proof}

\begin{lemma}
Each component of $F'$ is incompressible in $S^3$.
\end{lemma}

\begin{proof}
We specify the $\pm$-side of $F^{\pm}$ and $F'$ endowed from $F$ naturally.
Suppose without loss of generality that $F'$ is $+$-compressible, and let $E^+$ be a $+$-compressing disk for $F'$.
Then we can regard $E^+$ as a $\partial$-reducing disk for $E(F')$.
By applying our situation to \cite[Lemma 1.1.]{CG}, we may assume that $E^+$ intersects $H$ in a single loop which does not intersect $A_0$.
We deform $E^+$ by an isotopy so that $E^+\cap \mathcal{D}^+=\emptyset$ in $S^3$.
We take a complete meridian disk system $\mathcal{C}$ of $W$ which includes $\mathcal{D}^+$ and does not intersect $E^+$.
Put ${\mathcal{C}}^-={\mathcal{C}}-{\mathcal{D}}^+$.
Then we have $c(\rho (F;\partial E^+\cup \partial \mathcal{D}^+\cup \partial \mathcal{C}^-))$ $<$ $c(\rho (F;\partial \mathcal{D}^+\cup \partial \mathcal{C}^-))$ since $\partial E^+$ is essential in $F'$.
Suppose that $c(\rho (F;\partial E^+\cup \partial \mathcal{D}^+\cup \partial \mathcal{C}^-))$ $=$ $c(\rho (F;\partial E^+\cup \partial \mathcal{D}^+))$.
Then each component of $\partial \mathcal{D}^-$ bounds both $+$-compressing disk and $-$-compressing disk for $F$. 
Hence $F$ is reducible, but this contradicts Proposition 2.3.
Similarly, if $c(\rho (F;\partial E^+\cup \partial \mathcal{D}^+\cup \partial \mathcal{C}^-))$ $=$ $c(\rho (F;\partial \mathcal{C}^-))$, then we are done.
Hence we obtain a ${\pm}$-compressing disk system $E^+\cup \mathcal{D}^+$, $\mathcal{C}^-$ for $F$ which satisfies the conditions (1), (2) and have more minimal complexity than $\mathcal{D}^+\cup \mathcal{D}^-$.
This contradicts the property of $\mathcal{D}^+\cup \mathcal{D}^-$.

\end{proof}

By Lemma 7, we can deform the decomposing sphere $S$ by an isotopy so that $S$ intersects $F'$ in a single arc. Put $E(S)=S\cap E(F')$. Then $E(S)$ is a $\partial$-reducing disk for $E(F')$. Otherwise, at least one of $K_1$ and $K_2$ is trivial, and Theorem 1 clearly holds.
Let $V_0$ and $W_0$ be components of $V$ and $W$ respectively, where $V_0$ contains $A$ and $W_0$ is the next handlebody to $V_0$.
Put $H_0=V_0\cap W_0$.
Then $H_0$ gives a Heegaard splitting of $V_0\cup W_0$.
By Lemma 7, we can deform $E(S)$ by an isotopy so that $E(S)$ is contained in $V_0\cup W_0$.
By applying this situation to \cite[Lemma 1.1]{CG} or \cite[Theorem 1.3]{HS}, we may assume that $E(S)$ intersects $H_0$ in a single loop without moving $\partial E(S)$.
Moreover, there exist a complete meridian disk system ${\mathcal{E}}_0$ of $V_0$ such that ${\mathcal{E}}_0 \cap E(S)=\emptyset$ and ${\mathcal{E}}_0 \cap A_0=\emptyset$.
Thus $S$ intersects $F$ in a single arc, hence we have the conclusion.


\begin{thebibliography}{99}

\bibitem{CG}A.~J.~Casson and C.~McA.~Gordon, {\em Reducing Heegaard splittings}, Topology and its Appl. {\bf 27} (1987) 275-283.

\bibitem{H}J.~P.~Hempel, {\em 3-Manifolds}, Volume {\bf 86} of Ann. of Math. Stud. (Princeton Univ. Press, 1976).

\bibitem{HS}C.~Hayashi and K.~Shimokawa, {\em Thin position for 1-subminifold in 3-manifold}, preprint.

\bibitem{J}W.~H.~Jaco, {\em Lactures on Three-manifold Topology}, Volume {\bf 43} of CBMS Regional Conf. Ser. in Math. (American Math. Soc., 1980).

\bibitem{K}R.~Kirby, {\em problemlems in low-dimensional topology}, Part {\bf 2} of Geometric Topology (ed. W.~H.~Kazez), Studies in Adv. Math., (Amer. Math. Soc. Inter. Press, 1997).

\bibitem{KK}M.~Kobayashi and T.~Kobayashi, {\em On canonical genus and free genus of knot}, J. Knot Theory and Its Ramifi. {\bf 5} (1996) 77-85.

\bibitem{M}Y.~Moriah, {\em The free genus of knots,} Proc.~Amer.~Math.~Soc. {\bf 99} (1987) 373-379.

\bibitem{O}M.~Ozawa, {\em Synchronism of an incompressible non-free Seifert surface for a knot and an algebraically split closed incompressible surface in the knot complement}, to appear in Proc. Amer. Math. Soc.

\bibitem{S}H.~Schubert, {\em Knoten und Vollringe}, Acta Math. {\bf 90} (1953) 131-286.

\end{thebibliography}
\end{document}